\documentclass[a4paper,11pt]{amsart}
\usepackage{amsmath}
\usepackage{amssymb}
\usepackage{amsthm}
\usepackage{mathrsfs}
\usepackage{amsfonts}
\usepackage{verbatim}
\usepackage{geometry}
\usepackage{ulem}
\usepackage{indentfirst}
\usepackage{txfonts}
\def\l{\left}
\def\r{\right}
\def\f{\frac}

\def\ez{\epsilon}
\def\dz{\delta}

\def\et{\eta}
\def\lz{\lambda}
\def\Lz{\Lambda}
\def\az{\alpha}
\def\bz{\beta}

\def\tz{\theta}
\def\rz{\rho}
\def\Lfz{L^\infty}

\author{Lin Tang and Qian Zhang }
\title{\bf{Interior $C^{1,\az}$ regularity on the linearized Monge-Amp$\grave{e}$re equation with $\mathrm{VMO}$ type coefficients}\footnotetext{ \hspace{-0.65
cm} 2010 Mathematics Subject  Classification: 35J55.\\
The  research was supported  by the NNSF  (11271024)
and (11571289)  of China.}}
\date{}

\theoremstyle{plain}
\theoremstyle{plain}\newtheorem{thm}{Theorem}[section]
\theoremstyle{plain}\newtheorem{prop}{Proposition}[section]
\theoremstyle{plain}
\theoremstyle{plain}\newtheorem{lem}{Lemma}[section]
\theoremstyle{plain}\newtheorem{rem}{Remark}[section]
\newtheorem*{thm1}{Theorem 1}

\numberwithin{equation}{section}

\begin{document}
\maketitle
\noindent {\bf{Abstract}.}\quad
In this paper, we establish interior $C^{1,\az}$ estimates for solutions of the linearized Monge-Amp$\grave{e}$re equation
$$\mathcal{L}_{\phi}u:=\mathrm{tr}[\Phi D^2 u]=f,$$
where the density of the Monge-Amp$\grave{e}$re measure $g:=\mathrm{det}D^2\phi$ satisfies a $\mathrm{VMO}$-type condition and $\Phi:=(\mathrm{det}D^2\phi)(D^2\phi)^{-1}$ is the cofactor matrix of $D^2\phi$.

\bigskip

\section{Introduction}\label{s1}

This paper is concerned with interior regularity of solutions of the linearized Monge-Amp$\grave{e}$re equation
\begin{equation}\label{s1:1}
\mathcal{L}_\phi u:=\mathrm{tr}[\Phi D^2 u]=f,
\end{equation}
where $\phi$ is a solution of the Monge-Amp$\grave{e}$re equation
\begin{equation}\label{s1:2}
\mathrm{det}D^2\phi=g,\quad\quad \lz\le g\le\Lz\quad\mathrm{in}\;\Omega,\end{equation}
for some constants $0<\lz\le\Lz<\infty$.

The operator $\mathcal{L}_\phi$ appears in several contexts including affine geometry, complex geometry and fluid mechanics, see for example \cite{B,L,N1,NS,TW1}. In particular, the authors in \cite{TW1} resolved Chern's conjecture in affine geometry concerning affine maximal hypersurfaces in $\mathbb{R}^3$.

Concerning regularity of \eqref{s1:1} a fundamental result is the Harnack inequality for nonnegative solutions of $\mathcal{L}_\phi u=0$ established in \cite{CG}, which yields interior H$\ddot{o}$lder continuity of solutions of \eqref{s1:1}. By using this result and perturbation arguments, Guti$\acute{e}$rrez and Nguyen established in \cite{GN1} interior $C^{1,\alpha}$ estimates for solutions of \eqref{s1:1} when $g\in C(\Omega)$. Their main result has the form
\begin{equation}\label{s1:GN1}
\|u\|_{C^{1,\alpha'}(\Omega')}\leq C\{\|u\|_{\Lfz(\Omega)}+[f]^{n}_{\alpha,\Omega}\},
\end{equation}
for any $0<\alpha'<\alpha$. Here $\Omega'\Subset\Omega$ and $[f]^{n}_{\alpha,\Omega}$ is defined in Theorem 1 below.

In \cite{GN2}, interior $W^{2,p}$ estimates for solutions of \eqref{s1:1} were established for general $1<p<q$, $f\in L^q (q>n)$ and continuous density $g=\mathrm{det}D^2\phi$. The main result in this paper has the form
\begin{equation}\label{s1:GN2}
\|D^2 u\|_{L^{p}(\Omega')}\leq C\{\|u\|_{\Lfz(\Omega)}+\|f\|_{L^q(\Omega)}\}.
\end{equation}
By the imbedding theorem, when $1<n<p<q$ the above estimate holds if we replace the left hand side by $\|u\|_{C^{1,\gamma}(\Omega')}$ with $\gamma<1-n/p$. Since $[f]^{n}_{\alpha,\Omega}\leq C\|f\|_{L^q(\Omega)}$ for $1-\az-\f{n}{2q}\geq 0$, hence for $f\in L^q (q>n)$ the inequaity \eqref{s1:GN1} gives a better $C^{1,\gamma}$ estimate for $\gamma<1-n/p$ than \eqref{s1:GN2}.

On the other hand, in the case that $g$ is discontinuous, Huang \cite{H} proved interior $W^{2,p}$ estimates for solutions $\phi$ of \eqref{s1:2} where $g=\mathrm{det}D^2\phi$ belongs to a $\mathrm{VMO}$-type space $\mathrm{VMO}_\mathrm{loc}(\Omega,\phi)$ (see Section 2 for the definition). Using this result we recently established in \cite{LZ} global $W^{2,p}$ estimates for solutions of \eqref{s1:1} when $g\in\mathrm{VMO}_\mathrm{loc}(\Omega,\phi)$ defined in \cite{H}. Our result has a similar form to \eqref{s1:GN2} where the $L^p$ norm of $D^2u$ is estimated in terms of $\|f\|_{L^q(\Omega)}$. By the imbedding theorem again, interior $C^{1,\az}$ estimate when $g\in \mathrm{VMO}_\mathrm{loc}(\Omega,\phi)$ follows. But this $C^{1,\az}$ estimate is in terms of $\|f\|_{L^q(\Omega)}$ rather than $[f]^{n}_{\alpha,\Omega}$ in \eqref{s1:GN1}. Therefore, we are interested in establishing the interior $C^{1,\az}$ estimates for solutions of \eqref{s1:1} in terms of $[f]^{n}_{\alpha,\Omega}$ under the assumption that $g$ belongs to $\mathrm{VMO}_\mathrm{loc}(\Omega,\phi)$. Namely, we extend the result in \cite{GN1} from the case that $g\in C(\Omega)$ to the case that $g\in \mathrm{VMO}_\mathrm{loc}(\Omega,\phi)$. Our main result can be stated as follows.

\begin{thm1}\label{thm1}
Let $B_{\az_n}\subset\Omega\subset B_1$ be a normalized convex domain and $\phi\in C(\overline{\Omega})$ be a convex solution of \eqref{s1:2} with $\phi=0$ on $\partial\Omega$, where $g\in\mathrm{VMO}_{\mathrm{loc}}(\Omega,\phi)$. Assume that $u\in W^{2,n}_{\mathrm{loc}}(\Omega)$ is a solution of $\mathcal{L}_\phi u=f$ in $\Omega$ with $$[f]^{n}_{\alpha,\Omega}:=\sup_{S_\phi(x,r)\Subset\Omega}r^{\f{1-\az}{2}}\l(\f{1}{|S_\phi(x,r)|}\int_{S_\phi(x,r)}|f|^{n}dx\r)^{\f{1}{n}}<\infty$$
for some $0<\alpha<1$. Then for any $\alpha'\in(0,\alpha)$ and any $\Omega'\Subset \Omega$ we have
\begin{equation}
\|u\|_{C^{1,\az'}(\overline{\Omega'})}\leq C\{\|u\|_{\Lfz(\Omega)}+[f]^{n}_{\alpha,\Omega}\},
\end{equation}
where $C$ depends only on $n,\alpha,\alpha',\lambda,\Lambda,\mathrm{dist}(\Omega',\partial\Omega)$ and the $\mathrm{VMO}$-type property of $g$. \end{thm1}

The space $\mathrm{VMO}_{\mathrm{loc}}(\Omega,\phi)$ above is defined in Section 2.

We follows the perturbation arguments as in \cite{GN1}. The main lemma in our case is the stability of the cofactor matrix of $D^2\phi$ under a $\mathrm{VMO}$-type condition of $g=\mathrm{det}D^2\phi$. For this we use the interior $W^{2,p}$ estimates for solutions of \eqref{s1:2} in \cite{H}. We also need a result from \cite{H} which concerns the eccentricity of sections of \eqref{s1:2} under the $\mathrm{VMO}$-type condition of $g$.

The paper is organised as follows. In Section $2$, we establish the stability of cofactor matrix of $D^2\phi$ under a $\mathrm{VMO}$-type condition of $g=\mathrm{det}D^2\phi$. In Section $3$, we give an approximation lemma and investigate the eccentricity of sections of solutions of \eqref{s1:2} when $g$ is in $\mathrm{VMO}$-type spaces. In Section $4$, we prove the $C^{1,\az}$ estimate of solutions of \eqref{s1:1} at the minimum point of $\phi$, and finally, we give the complete proof of Theorem 1.

\section{Preliminary results and stability of cofactor matrices}\label{s2}

We first introduce some notation.

Let $\phi\in C(\overline{\Omega})$ be a solution of \eqref{s1:2}. A section of $\phi$ centered at $x_0\in\Omega$ with height $h$ is defined by
$$S_{\!\phi}(x_0,h):=\{x\in\Omega:\phi(x)<\phi(x_0)+\nabla\phi(x_0)\cdot(x-x_0)+h\}.$$
If $\phi=0$ on $\partial\Omega$, then for $0<\az<1$, we define
\begin{equation}\label{s2:GN24.41}
\Omega_\az:=\{x\in\Omega:\phi(x)<(1-\az)\min_\Omega\phi\}.
\end{equation}

Let $B_r(x_0)$ be the ball centered at $x_0\in\mathbb{R}^n$ with radius $r$ and denote for simplicity $B_r=B_r(0)$.

We always use the following assumption:\\
\noindent$\mathbf{(H)}$ $B_{a_1}\subset\Omega\subset B_{a_2}$ is a convex domain and $\phi\in C(\overline{\Omega})$ is a solution of \eqref{s1:2} with $\phi=0$ on $\partial\Omega$, where $0<a_1\leq a_2<\infty$.

Under the assumption $\mathbf{(H)}$ we often take $w$ to be the convex solution of
\begin{equation}\label{w}
\left\{
\begin{array}{rcl}
\mathrm{det}D^2 w=1&&{\mathrm{in}\;\Omega},\\
w=0&&{\mathrm{on}\;\partial\Omega}.
\end{array}\right.
\end{equation}

The following H$\ddot{o}$lder estimate for \eqref{s1:1} is from \cite{CG}.

\begin{lem}(See \cite[Lemma 2.5 and (2.2)]{GN1}.)\label{s2:lem GN12.5}
Assume that condition $\mathbf{(H)}$ holds. Let $u\in W^{2,n}_{\mathrm{loc}}(\Omega)$ be a solution of $\mathcal{L}_\phi u=f$ in $\Omega$. Let $\Omega'\Subset\Omega$ Then for any $x_0\in\Omega'$ and $h\le c$, we have
\begin{equation*}
|u(x)-u(y)|\leq C^*h^{-\bz}|x-y|^\bz\{\|u\|_{\Lfz(S_{\!\phi}(x_0,2h))}+(2h)^{\f{1}{2}}\|f\|_{L^n(S_{\!\phi}(x_0,2h))}\}
\quad\forall x,y\in S_{\!\phi}(x_0,h),
\end{equation*}
where $C^*,c>0, 0<\bz<1$ are constants depending only on $n,\lz,\Lz,a_1,a_2$ and $\mathrm{dist}(\Omega',\partial\Omega)$.
\end{lem}

The lemma below concerns classical $C^{1,1}$ interior estimate for uniformly elliptic equations.

\begin{lem}(See \cite[Theorem 2.7]{GN1}.)\label{s2:lem GN12.7}
Assume that $B_{a_1}\subset\Omega\subset B_{a_2}$ is a convex domain and $w$ is the solution of \eqref{w}. Then for any $\varphi\in C(\partial B_{\f{a_1}{2}})$ there exists a solution $h\in C^2(B_{\f{a_1}{2}})\cap C(\overline{B_{\f{a_1}{2}}})$ of $\mathcal{L}_w h=0$ in $B_{\f{a_1}{2}}$ and $h=\varphi$ on $\partial B_{\f{a_1}{2}}$ such that
\begin{eqnarray*}
\|h\|_{C^{1,1}(\overline{B_{\f{a_1}{4}}})}\leq c_e\|\varphi\|_{\Lfz(\partial B_{\f{a_1}{2}})},
\end{eqnarray*}
where $c_e>0$ is a constant depending only on $n,a_1,a_2$.
\end{lem}

The following Lemmas \ref{s2:lem H2.1}, \ref{s2:lem H3.1} and Theorem \ref{s2:thm HA(i)} were proved in \cite{H}.

\begin{lem}(See \cite[Lemma 2.1]{H}.)\label{s2:lem H2.1}
Assume that condition $\mathbf{(H)}$ holds. Then for any $\Omega'\Subset\Omega$, there exist positive constants $h_0, C$ and $q$ such that for $x_0\in\Omega'$, and $0<h\leq h_0$,
\begin{equation*}
B_{C^{-1}h}(x_0)\subset S_{\!\phi}(x_0,h)\subset B_{Ch^q}(x_0),
\end{equation*}
where $q=q(n,\lambda,\Lambda,a_1,a_2)$ and $h_0, C$ depends only on $n,\lambda,\Lambda, a_1, a_2$ and $\mathrm{dist}(\Omega',\partial\Omega)$.
\end{lem}

Assume that condition $\mathbf{(H)}$ holds. The space $\mathrm{VMO}_\mathrm{loc}(\Omega,\phi)$ is defined in \cite{H} as follows. Given a function $g\in L^1(\Omega)$, we say that $g\in \mathrm{VMO}_\mathrm{loc}(\Omega,\phi)$ if for any $\Omega'\Subset\Omega$,
$$Q_g(r,\Omega'):=\sup_{\substack{x_0\in\Omega',\\\mathrm{diam}(S_{\!\phi}(x_0,h))\leq r}}\mathrm{mosc}_{S_{\!\phi}(x_0,h)}\,g\rightarrow 0\quad\quad\mathrm{as}\;r\rightarrow 0.$$
Here the mean oscillation of $g$ over a measurable subset $A\subset\Omega$ is defined by
$$\mathrm{mosc}_A g:=\fint_A |g(x)-g_A|dx,$$
where $g_A=\fint_A g dx$ denotes the average of $g$ over $A$.

There are two simple facts about the definition above.

\begin{prop}\label{s2:prop VMO}
For any function $g^1$ such that $g:=(g^1)^n\in L^1(\Omega)$, the following hold:
\begin{enumerate}
\item[\rm(i)]
\begin{eqnarray*}
\left(\int_\Omega |g^1-(g^1)_\Omega|^n dx\right)^{\frac{1}{n}}\leq 2\left(\int_\Omega |g-g_\Omega| dx\right)^{\frac{1}{n}}.
\end{eqnarray*}
\item[\rm(ii)]
For any measurable subset $A\subset\Omega$, we have
\begin{eqnarray*}
\left(\int_\Omega |g^1-(g^1)_A|^n dx\right)^{\frac{1}{n}}\leq\left\{1+\left(\frac{|\Omega|}{|A|}\right)^{\frac{1}{n}}\right\}\left(\int_\Omega |g^1-(g^1)_\Omega|^n dx\right)^{\frac{1}{n}}.
\end{eqnarray*}
\end{enumerate}
\end{prop}

The maximum principle below is used to compare solutions $\phi$ of \eqref{s1:2} and $w$ of \eqref{w}, where $g=\mathrm{det}D^2\phi\in\mathrm{VMO}_\mathrm{loc}(\Omega,\phi)$ defined above.

\begin{lem}(See \cite[Lemma 3.1]{H}.)\label{s2:lem H3.1}
Assume $\Omega$ is a bounded convex domain in $\mathbb{R}^n$. Let $\phi$ and $w$ be the weak solutions to $\mathrm{det}D^2\phi=g_1^n\geq 0$ and $\mathrm{det}D^2 w=g_2^n\geq 0$ in $\Omega$, respectively. Assume that $g_1,g_2\in L^n(\Omega)$. Then
$$\max_{\overline{\Omega}}(\phi-w)\leq \max_{\partial\Omega}(\phi-w)+C_n\mathrm{diam}(\Omega)\left(\int_\Omega(g_2-g_1)^{+n}dx\right)^{1/n}.$$
\end{lem}

The theorem below gives $W^{2,p}$ estimates of solutions $\phi$ of \eqref{s1:2} under a $\mathrm{VMO}$-type condition of $g$. 

\begin{thm}(See \cite[Theorem A(i)]{H}.)\label{s2:thm HA(i)}
Assume the condition $\mathbf{(H)}$ holds. Let $0<\az<1, 1\leq p<\infty$. There exist constants $0<\ez<1$ and $C>0$ depending only on $n,\lz,\Lz,a_1,a_2, p$ and $\az$ such that if $\mathrm{mosc}_{S}g\leq\ez$ for any $S=S_{\!\phi}(x_0,h)\Subset\Omega$, then
$$\|D^2\phi\|_{L^p(\Omega_{\az})}\leq C.$$
\end{thm}

Next we establish stability of cofactor matrices of $D^2\phi$ under a $\mathrm{VMO}$-type condition of $g$.

\begin{lem}\label{s2:lem GN13.4}
Assume that $0<\lz\leq\Lz<\infty$ and $a_1,a_2>0$. Let $B_{a_1}\subset\Omega^k\subset B_{a_2}$ be a sequence of convex domain converging in the Hausdorff metric to a convex domain $B_{a_1}\subset\Omega\subset B_{a_2}$. For each $k\in\mathbb{N}$, let $\phi_k\in C(\overline{\Omega^k})$ be a convex function satisfying
$$\left\{
\begin{array}{rcl}
\mathrm{det}D^2\phi_k=g_k&&{\mathrm{in}\;\Omega^k},\\
\phi_k=0&&{\mathrm{on}\;\partial\Omega^k}.
\end{array}\right.$$
where $0<\lambda\leq g_k=(g^1_k)^n\leq\Lambda$ in $\Omega^k$, 
$$\mathrm{mosc}_{\Omega^k}\,g_k\leq \frac{1}{k}\quad\quad\mathrm{and}\quad\quad\sup_{S_{\!\phi_k}(x,h)\Subset\Omega^k}\,\mathrm{mosc}_{S_{\!\phi_k}(x,h)}g_k\leq\f{1}{k}.$$ 
Suppose that $\phi_k$ converges uniformly on compact subsets of $\Omega$ to a convex function $\phi\in C(\overline{\Omega})$ which is a solution of
$$\left\{
\begin{array}{rcl}
\mathrm{det}D^2\phi=1&&{\mathrm{in}\;\Omega},\\
\phi=0&&{\mathrm{on}\;\partial\Omega}.
\end{array}\right.$$
Then there exists a subsequence which we still denote by $\phi_{k}$ such that for any $1\leq p<\infty$,
\begin{equation*}
\lim_{k\rightarrow\infty}\|D^2\phi_{k}-(g^1_{k})_{B_{\f{a_1}{2}}}D^2\phi\|_{L^p(B_{\f{a_1}{2}})}=0,
\end{equation*}
and
\begin{equation*}
\lim_{k\rightarrow\infty}\|\Phi_{k}-(g^1_{k})_{B_{\f{a_1}{2}}}^{n-1}\Phi\|_{L^p(B_{\f{a_1}{2}})}=0,
\end{equation*}
where $\Phi_{k}$ and $\Phi$ are the cofactor matrices of $D^2\phi_{k}$ and $D^2\phi$ respectively.

\begin{proof}
First we note that since $\mathrm{dist}(B_{\f{a_1}{2}},\partial\Omega^k)\geq c(n,a_1,a_2)$, then $B_{\f{a_1}{2}}\subset\Omega^k_{\az}$ where $\az$ is a constant depending only on $n,\lz,\Lz,a_1,a_2$. For any $1\leq p<\infty$, let $\ez(p)=\ez(n,\lz,\Lz,a_1,a_2,p,\az)=\ez(n,\lz,\Lz,a_1,a_2,p)$ be the constant in Theorem \ref{s2:thm HA(i)}, then for any $k\geq k_{\ez(p)}:=\l[\f{1}{\ez(p)}\r]+1$ we have
\begin{equation}\label{s2:GN13.4*}
\sup_{S_{\!\phi_k}(x,h)\Subset\Omega^k}\mathrm{mosc}_{S_{\!\phi_k}(x,h)}\,g_k\leq\ez(p).
\end{equation}
Thus Theorem \ref{s2:thm HA(i)} implies that
\begin{equation}\label{s2:GN13.41}
\|D^2\phi_k\|_{L^p(B_{\f{a_1}{2}})}\leq\|D^2\phi_k\|_{L^p(\Omega^k_{\az})}\leq C(n,\lambda,\Lambda,a_1,a_2,p)\quad\forall k\geq k_{\ez(p)}.
\end{equation}

Let $\dz>0$ be an arbitrary small constant, and let $\Omega(\dz):=\{x\in\Omega:\mathrm{dist}(x,\partial\Omega)>\dz\}$. Then there exists $k_\dz\in\mathbb{N}$ such that for all $k\geq k_\dz$,
$$\mathrm{dist}(x,\partial\Omega^k)\leq 2\dz,\quad\quad\forall x\in\partial(\Omega(\dz)).$$
Then Aleksandrov's estimate (\cite[Theorem 1.4.2]{G}) implies that 
$$|\phi_k(x)-(g^1_k)_{B_{\f{a_1}{2}}}\phi(x)|\leq C(n,\lambda,\Lambda,a_1,a_2)\dz^{1/n}\quad\quad\forall x\in\partial(\Omega(\dz)).$$
By choosing $k_\dz$ even larger, we have $\Omega(\dz)\subset\Omega^k$ for $k\geq k_\dz$. It follows from Proposition \ref{s2:prop VMO} that,
\begin{eqnarray}\label{s2:GN13.45}
\left(\int_{\Omega(\dz)} |g^1_k-(g^1_k)_{B_{\f{a_1}{2}}}|^n dx\right)^{\frac{1}{n}}&\leq&\left(\int_{\Omega^k} |g^1_k-(g^1_k)_{B_{\f{a_1}{2}}}|^n dx\right)^{\frac{1}{n}}\leq C(n,a_1,a_2)\left(\int_{\Omega^k} |g^1_k-(g^1_k)_{\Omega^k}|^n dx\right)^{\frac{1}{n}}\nonumber\\
&\leq&C(n,a_1,a_2)\left(\int_{\Omega_k} |g_k-(g_k)_{\Omega_k}| dx\right)^{\frac{1}{n}}\leq \frac{C(n,a_1,a_2)}{k^{\f{1}{n}}}.
\end{eqnarray}

Using the above two estimates and applying Lemma \ref{s2:lem H3.1} with $\phi\rightsquigarrow\phi_k, w\rightsquigarrow (g^1_k)_{B_{\f{a_1}{2}}}\phi$, we get
\begin{eqnarray}\label{s2:GN13.42}
\max_{\overline{\Omega(\dz)}}|\phi_k-(g^1_k)_{B_{\f{a_1}{2}}}\phi|
&\leq& \max_{\partial(\Omega(\dz))}|\phi_k-(g^1_k)_{B_{\f{a_1}{2}}}\phi|
+C_n\mathrm{diam}(\Omega(\dz))\left(\int_{\Omega(\dz)}|g^1_k-(g^1_k)_{B_{\f{a_1}{2}}}|^{n}dx\right)^{\f{1}{n}}\nonumber\\
&\leq&C(n,\lambda,\Lambda,a_1,a_2)\l[\dz^{\f{1}{n}}+\f{1}{k^{\f{1}{n}}}\r],
\end{eqnarray}
for all $k\geq k_\dz$.

Using \eqref{s2:GN13.41}, \eqref{s2:GN13.45}, \eqref{s2:GN13.42} and similar arguments to the proof of \cite[Lemma 3.4]{GN1}, we can obtain the first conclusion of the lemma. For the second conclusion, we write

Write
$$\Phi_{k}-(g^1_{k})_{B_{\f{a_1}{2}}}^{n-1}\Phi=\l[1-\f{(g^1_{k})_{B_{\f{a_1}{2}}}^{n}}{\mathrm{det}D^2\phi_{k}}\r]\Phi_{k}
-\f{(g^1_{k})_{B_{\f{a_1}{2}}}^{n-1}}{\mathrm{det}D^2\phi_{k}}\Phi_{k}\l(D^2\phi_{k}-(g^1_{k})_{B_{\f{a_1}{2}}}D^2\phi\r)\Phi.$$
For any $1\leq q,r<\infty$, if $qr\leq n$ then by \eqref{s2:GN13.45} and H$\ddot{o}$lder inequality,
\begin{eqnarray*}
\l(\int_{B_{\f{a_1}{2}}}|g^1_{k}-(g^1_{k})_{B_{\f{a_1}{2}}}|^{qr}dx\r)^{\f{1}{qr}}&\leq& C(n,\lz,\Lz,a_1,a_2,q,r)\l(\int_{B_{\f{a_1}{2}}}|g^1_{k}-(g^1_{k})_{B_{\f{a_1}{2}}}|^{n}dx\r)^{\f{1}{n}}\\
&\leq&\f{C(n,\lz,\Lz,a_1,a_2,q,r)}{k^{\f{1}{n}}}.
\end{eqnarray*}
On the other hand, if $qr>n$ then
\begin{eqnarray*}
\l(\int_{B_{\f{a_1}{2}}}|g^1_{k}-(g^1_{k})_{B_{\f{a_1}{2}}}|^{qr}dx\r)^{\f{1}{qr}}&\leq&C(n,q,r,\lz,\Lz)\l(\int_{B_{\f{a_1}{2}}}|g^1_{k}-(g^1_{k})_{B_{\f{a_1}{2}}}|^{n}dx\r)^{\f{1}{qr}}\\
&\leq&\f{C(n,q,r,\lz,\Lz,a_1,a_2)}{k^{\f{1}{qr}}}.
\end{eqnarray*}
Note that 
\begin{eqnarray*}
|(g^1_{k})^{n}-(g^1_{k})_{B_{\f{a_1}{2}}}^{n}|\leq C(n,\lz,\Lz)|g^1_{k}-(g^1_{k})_{B_{\f{a_1}{2}}}|.
\end{eqnarray*}
Therefore, 
\begin{eqnarray*}
\l\|1-\f{(g^1_{k})_{B_{\f{a_1}{2}}}^{n}}{\mathrm{det}D^2\phi_{k}}\r\|_{L^{qr}(B_{\f{a_1}{2}})}\rightarrow 0,\quad \mathrm{as}\;k\rightarrow\infty.
\end{eqnarray*}
The rest of the proof is similar to that of \cite[Lemma 3.5]{GN1}, using \eqref{s2:GN13.41} and the first conclusion of the lemma.
\end{proof}
\end{lem}

\section{Main lemmas}\label{s3}

\subsection{A approximation lemma}\label{s3.1}
Next we compare solutions $v$ of \eqref{s1:1} and $h$ of $\mathcal{L}_w h=0$. The lemma below can be proved using similar arguments as in \cite[Lemma 4.1]{GN1}. The difference is that we estimate $\|v-h\|_{L^\infty}$ in terms of $\|\Phi-(g^{\f{1}{n}})^{n-1}_{B_1}W\|_{L^n}$ rather than $\|\Phi-W\|_{L^n}$.

\begin{lem}\label{s3:lem GN14.1}
Let $\rz^*:[0,\infty)\rightarrow[0,\infty)$ be a nondecreasing continuous function with $\lim_{\ez\rightarrow 0^+}\rz^*(\ez)\rightarrow 0$. Assume the condition $\mathbf{(H)}$ holds and $w\in C(\overline{\Omega})$ is the solution of \eqref{w}.
Suppose $v\in W^{2,n}_{\mathrm{loc}}(B_{\f{a_1}{2}})\cap C(\overline{B_{\f{a_1}{2}}})$ is a solution of $\mathcal{L}_\phi v=f$ in $B_{\f{a_1}{2}}$ with $|v|\leq 1$ in $B_{\f{a_1}{2}}$, and $h\in W^{2,n}_{\mathrm{loc}}(B_{\f{a_1}{2}})\cap C(\overline{B_{\f{a_1}{2}}})$ is a solution of
$$\left\{
\begin{array}{rcl}
\mathcal{L}_w h=0&&{\mathrm{in}\;B_{\f{a_1}{2}}},\\
h=v&&{\mathrm{on}\;\partial B_{\f{a_1}{2}}},
\end{array}\right.$$
Assume that $v$ and $h$ have $\rz^*$ as a modulus of continuity in $\overline{B_{\f{a_1}{2}}}$. Then for any $0<\tau<\f{a_1}{2}$, we have
\begin{eqnarray*}
\|v-h\|_{L^\infty(B_{\f{a_1}{2}-\tau})}
\leq C(n,\lz,\Lz,a_1,a_2)\l\{\rz^*\l(\|\Phi-(g^1)^{n-1}_{B_{\f{a_1}{2}}}W\|^{\f{1}{2}}_{L^n(B_{\f{a_1}{2}})}\r)+\|f\|_{L^n(B_{\f{a_1}{2}})}\r\}
\end{eqnarray*}
provided that $\|\Phi-(g^1)^{n-1}_{B_{\f{a_1}{2}}}W\|_{L^n(B_{\f{a_1}{2}})}\leq\tau^2$. Here $\Phi$ and $W$ are the cofactor matrices of $D^2\phi$ and $D^2 w$ respectively.
\end{lem}

Using Lemma \ref{s3:lem GN14.1}, the stability of the cofactor matrix Lemma \ref{s2:lem GN13.4} and arguing as in \cite[Lemma 4.2]{GN1}, we obtain the following approximation lemma when $g=\mathrm{det}\,D^2\phi$ satisfies a $\mathrm{VMO}$-type small oscillation. 

\begin{lem}\label{s3:lem GN14.2}
Let $\rz:[0,\infty)\rightarrow[0,\infty)$ be a nondecreasing continuous function with $\lim_{\ez\rightarrow 0^+}\rz(\ez)\rightarrow 0$. Given $K>0,\;\epsilon>0$. Assume the condition $\mathbf{(H)}$ holds and $w\in C(\overline{\Omega})$ is the solution of \eqref{w}. Let $\varphi\in C(\partial B_{\f{a_1}{2}})$ have $\rz$ as a modulus of continuity on $\partial B_{\f{a_1}{2}}$ and satisfy $\|\varphi\|_{L^\infty(\partial B_{\f{a_1}{2}})}\leq K$. Then there exists $\dz=\dz(\ez,n,\rz,\lz,\Lz,a_1,a_2,K)>0$ such that if 
$$\mathrm{mosc}_\Omega g\leq\dz\quad\quad\mathrm{and}\quad\quad\sup_{S_{\!\phi}(x,h)\Subset\Omega}\mathrm{mosc}_{S_{\!\phi}(x,h)}g\leq\dz,$$
and $f\in L^{n}(B_{\f{a_1}{2}})$ with $\|f\|_{L^{n}(B_{\f{a_1}{2}})}\leq\dz$, then any classical solutions $v, h$ of
$$\left\{
\begin{array}{rcl}
\mathcal{L}_\phi v=f&&{\mathrm{in}\;B_{\f{a_1}{2}}},\\
v=\varphi&&{\mathrm{on}\;\partial B_{\f{a_1}{2}}},
\end{array}\right.\quad\quad\mathrm{and}\quad\quad\left\{
\begin{array}{rcl}
\mathcal{L}_w h=0&&{\mathrm{in}\;B_{\f{a_1}{2}}},\\
h=v&&{\mathrm{on}\;\partial B_{\f{a_1}{2}}},
\end{array}\right.$$
satisfy
\begin{eqnarray*}
\|v-h\|_{\Lfz(B_{\f{a_1}{2}})}\leq\epsilon.
\end{eqnarray*}
\end{lem}

\subsection{Eccentricity of cross sections}\label{s3.2}

The two lemmas below concern the eccentricity of sections of \eqref{s1:2} if $g$ is in $\mathrm{VMO}$-type spaces. They are slight modifications of \cite[Lemmas 4.1, 4.2]{H}. Note that \cite[Lemma 4.1]{H} gives an affine transformation $Tx=A(x-z_0)$. But $z_0$ is not necessarily the minimum point of $\phi$. This is not convenient when we prove the $C^{1,\az}$ estimate in Theorem \ref{s4:thm GN14.5} in the next section. In the following two lemmas we replace $z_0$ in \cite[Lemma 4.1]{H} by the minimum point of $\phi$.

\begin{lem}\label{s3:lem GN13.2}
Assume the condition $\mathbf{(H)}$ holds, where $\mathrm{mosc}_\Omega g\leq\ez$. Then there exist $c_0, C_0>0$ depending only on $n,\lz,\Lz,a_1,a_2$ and a positive definite matrix $M=A^t A$ satisfying
\begin{equation*}
\mathrm{det}M=1,\quad 0<c_0 I\leq M\leq C_0 I,
\end{equation*}
such that for $0<\mu\leq c_0$ and $\ez^{\f{1}{n}}\leq c_0\mu^2$, we have
\begin{equation*}
B_{(1-\dz)\sqrt{\f{2}{(g^1)_\Omega}}}\subset\mu^{-\f{1}{2}}T S_{\!\phi}(x_0,\mu)\subset B_{(1+\dz)\sqrt{\f{2}{(g^1)_\Omega}}},
\end{equation*}
where $\dz=C_0(\mu^{\f{1}{2}}+\mu^{-1}\ez^{\f{1}{2n}})$, $x_0\in\Omega$ is the minimum point of $\phi$ and $Tx=A(x-x_0)$.

\begin{proof}
Let $w$ be the solution of \eqref{w}. Then from Lemma \ref{s2:lem H3.1} and Proposition \ref{s2:prop VMO}, we obtain
\begin{equation*}
\max_{\bar{\Omega}}|\phi-(g^1)_\Omega w|\leq C(n,\lz,\Lz,a_1,a_2)\ez^{\f{1}{n}}.
\end{equation*}
Using this and arguing as in \cite[Lemma 3.2]{GN1} we can obtain that there exist constants $C, c_0>0$ depending only on $n,\lz,\Lz,a_1,a_2$ such that if $0<\mu\leq c_0, 0<\gamma\leq\f{3\mu}{4}$ and $0<\ez^{\f{1}{n}}\leq c_0\mu^2$, then
\begin{eqnarray}
S_{\!w}(x_0,\mu-C\ez^{\f{1}{2n}})\subset S_{\!\f{\phi}{(g^1)_\Omega}}(x_0,\mu)\subset S_{\!w}(x_0,\mu+C\ez^{\f{1}{2n}}),\label{s3:GN13.21}\\
\partial S_{\!w}(x_0,\mu+\gamma)\subset N_{\f{C\gamma}{\sqrt{\mu}}}(\partial S_{\!w}(x_0,\mu)),\quad \partial S_{\!w}(x_0,\mu-\gamma)\subset N_{\f{C\gamma}{\sqrt{\mu}}}(S_{\!w}(x_0,\mu)),\label{s3:GN13.22}\\
B_{C\sqrt{\mu}}(x_0)\subset S_{\!w}(x_0,\mu)\subset B_{C\sqrt{\mu}}(x_0),\label{s3:GN13.23}\\
\partial S_{\!w}(x_0,\mu)\subset N_{C\mu}(\partial\sqrt{\mu}E),\label{s3:GN13.24}
\end{eqnarray}
where $E:=\l\{x:\f{1}{2}\langle D^2 w(x_0)(x-x_0),x-x_0\rangle\leq 1\r\}$.

From \eqref{s3:GN13.21}, \eqref{s3:GN13.22} and \eqref{s3:GN13.24}, we obtain
$$\partial S_{\!\f{\phi}{(g^1)_\Omega}}(x_0,\mu)\subset N_{C\f{\ez^{\f{1}{2n}}}{\sqrt{\mu}}}(\partial S_{\!w}(x_0,\mu))\subset N_{C\l(\mu+\f{\ez^{\f{1}{2n}}}{\sqrt{\mu}}\r)}(\partial\sqrt{\mu}E).$$
Since $S_{\!\phi}(x_0,\mu)=S_{\!\f{\phi}{(g^1)_\Omega}}\l(x_0,\f{\mu}{(g^1)_\Omega}\r)$, we obtain
\begin{eqnarray}\label{s3:GN13.25}
\partial S_{\!\phi}(x_0,\mu)\subset N_{C\l(\mu+\f{\ez^{\f{1}{2n}}}{\sqrt{\mu}}\r)}\l(\partial\sqrt{\f{\mu}{(g^1)_\Omega}}E\r)
\end{eqnarray}
for any $0<\mu\leq c_0$ and $\ez^{\f{1}{n}}\leq c_0\mu^2$.

Write $D^2 w(x_0)=A^t A$ for some positive definite matrix $A>0$ and $M=D^2 w(x_0)$, and then the conclusion follows.
\end{proof}
\end{lem}

\begin{lem}\label{s3:lem GN13.3}
Assume that $\lz\leq b^n\leq\Lz$. Let $B_{(1-\dz)\sqrt{\f{2}{b}}}\subset\Omega\subset B_{(1+\dz)\sqrt{\f{2}{b}}}$ be a convex domain with $\dz>0$ small and $\phi\in C(\overline{\Omega})$ be a solution of \eqref{s1:2} with $\phi=0$ on $\partial\Omega$, where $\mathrm{mosc}_\Omega g\leq\ez$. Then there exist $c_0, C_0>0$ depending only on $n,\lz,\Lz$ and a positive definite matrix $M=A^t A$ satisfying
\begin{equation*}
\mathrm{det}M=1,\quad (1-C_0\dz)I\leq M\leq (1+C_0\dz)I,
\end{equation*}
such that for $0<\mu\leq c_0$ and $\ez^{\f{1}{n}}\leq c_0\mu^2$, we have
\begin{equation*}
B_{(1-\dz_1)\sqrt{\f{2}{(g^1)_\Omega}}}(0)\subset\mu^{-\f{1}{2}}T S_{\!\phi}(x_0,\mu)\subset B_{(1+\dz_1)\sqrt{\f{2}{(g^1)_\Omega}}}(0),
\end{equation*}
where $\dz_1=C_0(\dz\mu^{\f{1}{2}}+\mu^{-1}\ez^{\f{1}{2n}})$, $x_0\in\Omega$ is the minimum point of $\phi$ and $Tx=A(x-x_0)$.

\begin{proof}
As in Lemma \ref{s3:lem GN13.2}, let $w$ be the solution of \eqref{w}. Then \eqref{s3:GN13.21}, \eqref{s3:GN13.22}, \eqref{s3:GN13.23} still hold. Let 
$$E:=\l\{x:\f{1}{2}\langle D^2 w(x_0)(x-x_0),x-x_0\rangle\leq 1\r\}$$ 
be as in the proof of Lemma \ref{s3:lem GN13.2}. Then instead of \eqref{s3:GN13.24}, we have
\begin{eqnarray}\label{s3:GN13.31}
\partial S_{\!w}(x_0,\mu)\subset N_{C\dz\mu}(\partial\sqrt{\mu}E).
\end{eqnarray}
Indeed, we argue as in the proof of \cite[(3.16)]{GN1} and find that in order to prove \eqref{s3:GN13.31}, we only need to prove that for any $\xi\in(1+C\dz\sqrt{\mu})\sqrt{\mu}E$, we have
\begin{eqnarray}\label{s3:GN13.32}
|D^3 w(\xi)|\leq C(n,\lz,\Lz)\dz.
\end{eqnarray}
For this, we note that $(1+C\dz\sqrt{\mu})\sqrt{\mu}E\subset\Omega'\Subset\Omega$ for some $\Omega'$. Then from the proof of \cite[Lemma 4.2]{H}, we obtain
\begin{eqnarray*}
\|w-P\|_{C^3(\overline{\Omega'})}\leq C(n,\lz,\Lz)\dz,
\end{eqnarray*}
where $P(x)=\f{1}{2}|x|^2-\f{1}{b}$. Thus, the estimate \eqref{s3:GN13.32} holds and therefore \eqref{s3:GN13.31} is true. Moreover, the last estimate implies that
\begin{eqnarray}\label{s3:GN13.33}
|D^2 w(x_0)-I|\leq C(n,\lz,\Lz)\dz.
\end{eqnarray}

Similar to Lemma \ref{s3:lem GN13.2}, \eqref{s3:GN13.21}, \eqref{s3:GN13.22} and \eqref{s3:GN13.31} imply that
$$\partial S_{\!\f{\phi}{(g^1)_\Omega}}(x_0,\mu)\subset N_{C\f{\ez^{\f{1}{2n}}}{\sqrt{\mu}}}(\partial S_{\!w}(x_0,\mu))\subset N_{C\l(\dz\mu+\f{\ez^{\f{1}{2n}}}{\sqrt{\mu}}\r)}(\partial\sqrt{\mu}E)$$
or
\begin{eqnarray}\label{s3:GN13.34}
\partial S_{\!\phi}(x_0,\mu)\subset N_{C\l(\dz\mu+\f{\ez^{\f{1}{2n}}}{\sqrt{\mu}}\r)}\l(\partial\sqrt{\f{\mu}{(g^1)_\Omega}}E\r)
\end{eqnarray}
for any $0<\mu\le c_0$ and $\ez^{\f{1}{n}}\leq c_0\mu^2$.

Write $D^2 w(x_0)=A^t A$ for some positive definite matrix $A>0$ and $M=D^2 w(x_0)$, then \eqref{s3:GN13.33} gives
$$(1-C\dz)I\leq M\leq (1+C\dz)I,$$
and the conclusion follows.
\end{proof}
\end{lem}

\begin{rem}\label{s3:rem GN13.3}
Under the assumptions in Lemmas \ref{s3:lem GN13.2} and \ref{s3:lem GN13.3}, it follows from \eqref{s3:GN13.21} and \eqref{s3:GN13.23} that for any $0<\mu\leq c_0$ and $\ez^{\f{1}{n}}\leq c_0\mu^2$, we have
\begin{equation}\label{s3:GN13.3*}
S_{\!\phi}(x_0,\mu)\subset B_{C\sqrt{\mu+C\ez^{\f{1}{2n}}}}(x_0),
\end{equation}
where $C, c_0$ depend only on $n,\lz,\Lz$ (Under the assumptions in Lemma \ref{s3:lem GN13.2} these constants also depend on $a_1,a_2$).
\end{rem}

\section{Interior $C^{1,\az}$ estimate for linearized equation}\label{s4}

\subsection{Estimate at the minimum point of the convex function}
In this subsection, we prove $C^{1,\az}$ estimate of \eqref{s1:1} at the minimum point of $\phi$ under a $\mathrm{VMO}$-type condition of $\mathrm{det}\,D^2\phi$.

\begin{thm}\label{s4:thm GN14.5}
Assume that $0<\alpha'<\alpha<1, r_0, C_1>0$ and $0<\lz\leq\Lz<\infty$. Assume $B_{\az_n}\subset\Omega\subset B_{1}$ is a normalized convex domain and $\phi\in C(\overline{\Omega})$ is a convex solution of \eqref{s1:2} with $\phi=0$ on $\partial\Omega$, where 
$$\mathrm{mosc}_\Omega g\leq\tz\quad\quad\mathrm{and}\quad\quad\sup_{S_{\!\phi}(x,h)\Subset\Omega}\mathrm{mosc}_{S_{\!\phi}(x,h)}g\leq\tz.$$
Let $u\in W^{2,n}_{\mathrm{loc}}(\Omega)$ be a solution of $\mathcal{L}_\phi u=f$ in $\Omega$ with
\begin{equation*}
\l(\f{1}{|S_{\!r}(\phi)|}\int_{S_{\!r}(\phi)}|f|^{n}dx\r)^{\f{1}{n}}\leq C_1 r^{\f{\alpha-1}{2}}\quad\quad \mathrm{for\;all}\;S_{\!r}(\phi)=S_{\!\phi}(x_0,r)\Subset\Omega, r\leq r_0,
\end{equation*}
where $x_0$ is the minimum point of $\phi$, then $u$ is $C^{1,\az'}$ at $x_0$, more precisely, there is an affine function $l(x)$ such that
\begin{equation*}
r^{-(1+\alpha')}\|u-l\|_{\Lfz(B_r(x_0))}+|Dl|\leq C\{|u\|_{\Lfz(\Omega)}+C_1\}\quad\quad\forall r\leq \mu^*,
\end{equation*}
where $\tz\in(0,1), C>0,\;\mu^*>0$ depend only on $n,\lz,\Lz,\alpha,\alpha',r_0$.

\begin{proof}
We can assume that $\|u\|_{\Lfz(\Omega)}\leq 1$ and
\begin{equation}\label{s4:GN14.51}
\l(\f{1}{|S_{\!r}(\phi)|}\int_{S_{\!r}(\phi)}|f|^{n}dx\r)^{\f{1}{n}}\leq \tz r^{\f{\alpha-1}{2}}\quad\quad \mathrm{for\;all}\;S_{\!r}(\phi)=S_{\!\phi}(x_0,r)\Subset\Omega, r\leq r_0.
\end{equation}
And we only need to prove that
\begin{equation}\label{s4:GN14.52}
r^{-(1+\alpha')}\|u-l\|_{\Lfz(B_r(x_0))}+|Dl|\leq C\quad\quad\forall r\leq \mu^*,
\end{equation}
where $C>0,\;\mu^*>0$ depend only on $n,\lz,\Lz,\alpha,\alpha',r_0$.

Define $a_1:=\f{1}{2}\sqrt{\f{2}{\Lz^{\f{1}{n}}}}$ and $a_2:=2\sqrt{\f{2}{\lz^{\f{1}{n}}}}$. Then Lemma \ref{s2:lem GN12.5} gives constants $C^*,\bz>0$, Lemma \ref{s2:lem GN12.7} gives $c_e>0$, Lemma \ref{s3:lem GN13.2} and \ref{s3:lem GN13.3} give $c_0,C_0>0$. All these constants depend only on $n,\lz,\Lz$. Applying Lemma \ref{s3:lem GN13.2}, \ref{s3:lem GN13.3} and \ref{s3:lem GN14.2} and similar arguments to the proof of \cite[Theorem 4.5]{GN1}, we can prove that there exist $0<\mu<1$ depending only on $n,\lz,\Lz,\alpha,r_0$, a sequence of positive matrices $A_k$ with $\mathrm{det}A_k=1$, a sequence $b_k>0$ and
a sequence of affine functions $l_k(x)=a_k+B_k\cdot(x-x_0)$ satisfying for $k\geq 1$,
\begin{equation*}
\|A_{k-1}A_k^{-1}\|\leq\f{1}{\sqrt{c_0}},\quad\|A_k\|\leq\sqrt{C_0(1+C_0\dz_0)(1+C_0\dz_1)\cdots(1+C_0\dz_{k-1})};\tag{1}
\end{equation*}

\begin{equation*}
B_{a_1}\subset B_{(1-\dz_k)\sqrt{\f{2}{b_k}}}\subset\mu^{-\f{k}{2}}A_k(S_{\!\mu^k}(\phi)-x_0)\subset B_{(1+\dz_k)\sqrt{\f{2}{b_k}}}\subset B_{a_2},\quad\quad\lz\leq b_k^n\leq\Lz;\tag{2}
\end{equation*}

\begin{equation*}
\|u-l_{k-1}\|_{\Lfz(S_{\!\mu^k}(\phi))}\leq\mu^{\f{1+\alpha}{2}(k-1)};\tag{3}
\end{equation*}

\begin{equation*}
|a_k-a_{k-1}|+\mu^{\f{k}{2}}|(A_k^{-1})^t (B_k-B_{k-1})|\leq 2c_e\mu^{\f{1+\alpha}{2}(k-1)};\tag{4}
\end{equation*}

\begin{equation*}
\f{|(u-l_{k-1})(\mu^{\f{k}{2}}A_k^{-1}x+x_0)-(u-l_{k-1})(\mu^{\f{k}{2}}A_k^{-1}y+x_0)|}{\mu^{\f{1+\alpha}{2}(k-1)}}
\leq 2 C^*(\sqrt{c_1\mu})^{-\bz}|x-y|^\bz,\tag{5}
\end{equation*}
for any $x,y\in\mu^{-\f{k}{2}}A_k(S_{\!\mu^k}(\phi)-x_0)$,

where $A_0:=I,\quad l_0:=0,\quad \dz_0:=0,\quad \dz_1:=C_0(\mu^{\f{1}{2}}+\mu^{-1}\tz^{\f{1}{2n}})<\f{1}{2},$\\

$\dz_k:=C_0(\dz_{k-1}\mu^{\f{1}{2}}+\mu^{-1}\tz^{\f{1}{2n}}),\quad \dz_k<\dz_{k-1}\quad \mathrm{for}\; k\geq 2.$\\

The rest of the proof is the same as Part $4$ (proof of (4.35)) in the proof of \cite[Theorem 4.5]{GN1}.
\end{proof}
\end{thm}

\subsection{Proof of Theorem 1}

By Lemma \ref{s2:lem H2.1}, for any $\Omega'\Subset\Omega$, there exist positive constants $h_0, C$ and $q$ depending only on $n,\lambda,\Lambda$ and $\mathrm{dist}(\Omega',\partial\Omega)$ such that for any $x_0\in\Omega'$, we have
\begin{equation}\label{s4:GN14.73}
B_{C^{-1}h_0}(x_0)\subset S_\phi(x_0,h_0)\subset B_{Ch_0^q}(x_0).
\end{equation}
Choose $h_0$ smaller and we can assume $S_\phi(x_0,h_0)\subset B_{C(h_0)^q}(x_0) \subset\Omega''\Subset\Omega$. Since $g\in \mathrm{VMO}_\mathrm{loc}(\Omega,\phi)$, we have
$$\et_g(r,\Omega''):=\sup_{\substack{S_{\!\phi}(x,h)\subset\Omega'',\\\mathrm{diam}(S_{\!\phi}(x,h))\leq r}}\mathrm{mosc}_{S_{\!\phi}(x,h)}\,g\rightarrow 0,\quad r\rightarrow 0.$$
Let $\tz=\tz(n,\az,\az',\lz,\Lz,\mathrm{dist}(\Omega',\partial\Omega))>0$ be the constant in Theorem \ref{s4:thm GN14.5}, then there exists $0<r_1<1$ such that $\et_g(r_1,\Omega'')<\tz$. Take $h_0$ smaller such that $\mathrm{diam}(B_{Ch_0^q}(x_0))\leq r_1$, then for any $S_{\!\phi}(x,h)\subset S_{\!\phi}(x_0,h_0)$, we have $S_{\!\phi}(x,h)\subset\Omega''$ and $\mathrm{diam}(S_{\!\phi}(x,h))\leq r_1$, thus,
\begin{equation}\label{s4:GN14.74}
\mathrm{mosc}_{S_{\!\phi}(x,h)}\,g\leq\et_g(r_1,\Omega'')\le\tz.
\end{equation}
Fix such $h_0$ in the rest of the proof. Note that $h_0$ depends only on $n,\lambda,\Lambda, \mathrm{dist}(\Omega',\partial\Omega),r_1$. Thus it depends only on $n,\lambda,\Lambda, \mathrm{dist}(\Omega',\partial\Omega),\az,\az'
$ and the $\mathrm{VMO}$-type property of $g$.

Let $T$ be an affine map such that
\begin{equation*}
B_{\az_n}\subset T(S_{\!\phi}(x_0,h_0))\subset B_1.
\end{equation*}
By \eqref{s4:GN14.73} we have
\begin{equation}\label{s4:GN14.75}
\|T\|\leq C h^{-1}_0,\quad \|T^{-1}\|\leq C h_0^q,
\end{equation}
where $C=C(n,\lambda,\Lambda, \mathrm{dist}(\Omega',\partial\Omega))>0$. Let $\kappa_0:=|\mathrm{det}\,A|^{\f{2}{n}}$, then we have $\kappa_0 h_0\geq r_0$ for some constant $r_0$ depending only on $n,\lz,\Lz$.

For $y\in \tilde{\Omega}:=T(S_{\!\phi}(x_0,h_0))$, define
$$\tilde{\phi}(y)=\kappa_0[(\phi-l_{x_0})(T^{-1}y)-h_0]\quad\quad\mathrm{and}\quad\quad v(y)=\kappa_0^{\f{1+\alpha}{2}}u(T^{-1}y),$$
where $l_{x_0}(x)$ is the supporting function of $\phi$ at $x_0$. Then,
$$\mathrm{det}D^2\tilde{\phi}(y)=\tilde{g}(y)=(\tilde{g}^1(y))^n,\quad\quad\lz\leq\tilde{g}(y)=g(T^{-1}y)\leq\Lz\quad \mathrm{in}\;\tilde{\Omega}$$
and by \eqref{s4:GN14.74}, we have
$$\mathrm{mosc}_{\tilde{\Omega}}\,\tilde{g}=\mathrm{mosc}_{S_{\!\phi}(x_0,h_0)}\,g\leq\tz$$
and
$$\sup_{S_{\!\tilde{\phi}}(y,h)\Subset \tilde{\Omega}}\mathrm{mosc}_{S_{\!\tilde{\phi}}(y,h)}\,\tilde{g}\leq\sup_{S_{\!\phi}(T^{-1}y,\,h\kappa_0^{-1})\subset S_{\!\phi}(x_0,h_0)}\mathrm{mosc}_{S_{\!\phi}(T^{-1}y,\,h\kappa_0^{-1})}\,g\leq\tz.$$
Applying Theorem \ref{s4:thm GN14.5} to $v$ and arguing as in the proof of \cite[Theorem 4.7]{GN1}, we obtain the conclusion of Theorem 1.

LMAM, School of Mathematical Sciences,
 Peking University, Beijing, 100871,
 P. R. China

 Lin Tang,\quad
 E-mail address:  tanglin@math.pku.edu.cn

Qian Zhang,\quad
E-mail address: 1401110018@pku.edu.cn

\begin{thebibliography}{99}

\bibitem{B}  \label{B} Brenier, Y.,
\textit{Polar factorization and monotone rearrangement of vector valued functions,}
Comm. Pure Appl. Math. $\mathbf{44}$ (1991), 375-417.



\begin{comment}
\bibitem{C}  \label{C} Caffarelli, L.A.,
\textit{Interior $W^{2,p}$ estimates for solutions of the Monge-Amp$\grave{e}$re Equations,}
Ann. Math. $\mathbf{131}$ (1990), 135-150.

\bibitem{C2}  \label{C2} Caffarelli, L.A.,
\textit{A localization property of viscosity solutions to the Monge-Amp$\grave{e}$re Equation,}
Ann. Math. $\mathbf{131}$ (1990), 129-134.

\bibitem{C4}  \label{C4} Caffarelli, L.A.,
\textit{Some regularity properties of solutions of Monge-Amp$\grave{e}$re Equation,}
Comm. Pure Appl. Math. $\mathbf{44}$ (1991), 965-969.

\bibitem{CC}  \label{CC} Caffarelli, L.A., Cabr$\acute{e}$, X.,
\textit{Fully Nonlinear Elliptic equations,}
 Amer. Math. Soc. Colloq. Publ.,$\mathbf{43}$, American Mathematical Society, Procidence, RI, (1995).
\end{comment}




\bibitem{CG}  \label{CG} Caffarelli, L.A., Guti$\acute{e}$rrez, C.E.,
\textit{Properities of the solutions of the linearized Monge-Amp$\grave{e}$re equations,}
 Amer. J. Math. $\mathbf{119}$ (1997), 423-465.




\bibitem{G}  \label{G} Guti$\acute{e}$rrez, C.E.,
\textit{The Monge-Amp$\grave{e}$re Equation,}
Birkh$\ddot{a}$user, Boston, MA, 2001.




\bibitem{GN1}  \label{GN1} Guti$\acute{e}$rrez, C.E., Nguyen, T.V.,
\textit{Interior gradient estimates for solutions to the linearized Monge-Amp$\grave{e}$re equation,}
 Adv. Math. $\mathbf{228}$ (4) (2011), 2034-2070.

\bibitem{GN2}  \label{GN2} Guti$\acute{e}$rrez, C.E., Nguyen, T.V.,
\textit{Interior second derivative estimates for solutions to the linearized Monge-Amp$\grave{e}$re equation,}
Trans. Amer. Math. Soc. $\mathbf{367}$ (8) (2015), 4537-4568.

\bibitem{GT}  \label{GT} Guti$\acute{e}$rrez, C.E., Tournier, F.,
\textit{$W^{2,p}$-estimates for the linearized Monge-Amp$\grave{e}$re equation,}
Trans. Amer. Math. Soc. $\mathbf{358}$ (11) (2006), 4843-4872.

\bibitem{H}  \label{H} Huang, Q.,
\textit{Sharp regularity results on second derivatives of solutions to the Monge-Amp$\grave{e}$re Equation with VMO type data,}
Commun. Pure. Appl. Math. $\mathbf{62}$ (5) (2009), 677-705.

\bibitem{L}  \label{L} Loeper, G.,
\textit{A fully nonlinear version of the incompressible Euler equations: the semigeostrophic system,}
SIAM J. Math. Anal. $\mathbf{38}$ (3) (2006), 795-823.

\bibitem{N1}  \label{N1} Le, N.Q.,
\textit{$W^{4,p}$ solution to the second boundary value problem of the priscribed affine mean curture and Abreu's Equations,}
arXiv:1506.06727 [math.AP].

\bibitem{NS}  \label{NS} Le, N.Q., Savin, O.,
\textit{Some minimization problems in the class of convex functions with priscribed determinant,}
Anal. PDE. $\mathbf{6}$ (5)(2013), 1025-1050.



\begin{comment}
\bibitem{NT}  \label{NT} Le, N.Q., Nguyen, T.V.,
\textit{Global $W^{2,p}$ estimates for solutions to the linearized Monge-Amp$\grave{e}$re Equations,}
Ann. Math. $\mathbf{358}$ (3-4)(2014), 629-700.

\bibitem{S1}  \label{S1} Le, N.Q., Savin, O.,
\textit{Boundary regularity for solutions to the linearized Monge-Amp$\grave{e}$re Equations,}
Arch. Ration. Mech. Anal. $\mathbf{210}$ (3)(2013), 813-836.

\bibitem{S2}  \label{S2} Savin, O.,
\textit{Pointwise $C^{2,\alpha}$ estimates at the boundary for the Monge-Amp$\grave{e}$re Equation,}
J. Amer. Math. Soc. $\mathbf{26}$ (1)(2013), 63-99.

\bibitem{S3}  \label{S3} Savin, O.,
\textit{Global $W^{2,p}$ estimates for the Monge-Amp$\grave{e}$re Equations,}
Proc. Am. Math. Soc. $\mathbf{141}$ (10)(2013), 3573-3578.

\bibitem{TW}  \label{TW} Trudinger, N.S., Wang, X.J.,
\textit{Boundary regularity for Monge-Amp$\grave{e}$re and affine maximal surface equations,}
Ann. Math. $\mathbf{167}$ (3)(2008), 993-1028.
\end{comment}



\bibitem{TW1}  \label{TW1} Trudinger, N.S., Wang, X.J.,
\textit{The Bernstein problem for affine maximal hypersurfaces,} Invent. Math, 140 (2000), 399-422.

\bibitem{LZ}  \label{LZ} Tang, L., Zhang, Q.,
\textit{Global $W^{2,p}$ regularity on the linearized Monge-Amp$\grave{e}$re equation with $\mathrm{VMO}$ type coefficients,}
submitted.
\end{thebibliography}
\end{document}